\documentclass[a4paper,10pt]{amsart}
\usepackage{epsf}  
\usepackage{amsfonts, amssymb, amsthm, amsmath}  
\usepackage{hyperref}
\usepackage{tikz}
\usetikzlibrary{matrix,arrows,decorations.pathmorphing}
\usepackage{tikz-cd}
\usepackage{graphicx}
\usepackage{psfrag}

\numberwithin{thm}{section}

\providecommand{\totl}[1]{\ensuremath{\lceil #1\rceil }}
\providecommand{\totb}[1]{\ensuremath{\underline{ #1}}}

\newcommand{\ex}{\bold}

\newcommand{\dbar}{\bar{\partial}}

\author{Brett Parker   }
\email{brettdparker@gmail.com}

\title[]{On the value of thinking tropically to understand Ionel's GW invariants relative normal-crossing divisors}

\thanks{This research was supported by the ARC grant DP1093094.}
\begin{document}
\maketitle

\begin{abstract} Ionel's GW invariants relative normal-crossing divisors appear different  from  Gromov--Witten invariants defined using log schemes or exploded manifolds. Appearances are, in this case, deceiving.  I sketch the relationship between Ionel's invariants and their exploded cousins using the  example of the moduli space of lines in the complex projective plane relative  two coordinate lines. Even in this simplest of examples,   13 different types of curves appear in Ionel's compactified moduli space, but these different types of curves can be understood in a unified and intuitive fashion using tropical curves. 

\end{abstract}

There are three separate approaches to defining Gromov--Witten invariants relative  normal-crossing divisors: log geometry\footnote{ Log Gromov--Witten invariants are defined by  Gross and Sibert in \cite{GSlogGW} and Abramovich and Chen in \cite{Chen,Chen2,acgw,acev}. Earlier work by Kim studied  log Gromov--Witten invariants in a less general situation \cite{kim}. In \cite{AMW}, Abramovich, Marcus, and Wise prove that, in the case of a smooth divisor, these log Gromov--Witten invariants give equivalent invariants to the relative Gromov--Witten invariants defined by Li in \cite{Li}, so log Gromov--Witten invariants are closely related to the invariants defined by Li and Ruan in \cite{ruan}, and Ionel and Parker \cite{IP1}.}, Ionel's GW invariants relative normal-crossing divisors, and my approach using exploded manifolds\footnote{I define Gromov--Witten invariants of exploded manifolds  in the series of papers \cite{iec,cem,reg,dre,uts,evc,vfc}, and prove the associated gluing formula  in \cite{gfgw}. For a brief introduction to exploded manifolds and the tropical gluing formula, see \cite{scgp}, and for examples of the tropical gluing formula in action, see \cite{3d, tec, tpgf}.   This gluing theorem is also described informally in a talk with lots of pictures --- found  \href{http://prezi.com/6o-kwhznakpf/tropical-gluing-formulas-for-Gromov--witten-invariants/}{at this link}, or by searching the internet for `Brett Parker tropical gluing prezi'.}. In the algebraic setting, the approach using log geometry and exploded manifolds is related by a base change\footnote{See \cite{elc} for an explanation of the relationship between log schemes and exploded manifolds. This close relationship was accidental, and forced by the geometry of holomorphic curves in normal-crossing degenerations.}, and, technicalities aside, the two approaches give equivalent invariants. What is less obvious is the relationship between Ionel's GW invariants relative normal-crossing divisors and their log or  exploded analogues. 

Ionel's definition of GW invariants relative normal-crossing divisors relies on the following observations:
\begin{itemize}
\item Given an appropriate generic perturbation of the $\dbar$ equation, there is a compactification of the moduli space of (perturbed) holomorphic curves not contained in the divisor, and this compactification adds strata that are  smooth manifolds with codimension at least $2$.
\item The appropriate evaluation map is continuous, and is also smooth when restricted to each strata.
\end{itemize}
 In fact, the same perturbation can achieve transversality in Ionel's setting and the exploded manifold setup. So,  these two approaches define roughly the same invariants, despite differences between the two compactifications.
The compactification defined by Ionel is a kind of blow up of the compactification achieved using exploded manifolds, but all evaluation maps used by Ionel actually factor through the smaller `exploded' compactification.  In this note, we  compare these two compactifications using  the easiest instructive example. In doing so, I hope to emphasize the value of thinking tropically in order to understand Ionel's compactification.

%
%
%
%
%
%
%

\

Consider $\mathbb CP^{2}$ with the normal-crossing divisor $D_{1}\cup D_{2}$,  where $D_{i}$ is $\{z_{i}=0\}$. This structure is not quite enough to define $D_{1}\cup D_{2}$ as a normal-crossing divisor in the sense of Ionel's Definition 1.3 of \cite{IonelGW}\footnote{All references to \cite{IonelGW} are to the published version, and version 3 on the arXiv; this will be important later, because one of the improvements made between version 1 and 3 is that in version 3,  curves in the compactified moduli space must satisfy an extra condition equivalent to the existence of a certain tropical curve.},
 however, if we make the obvious identifications of $\mathbb CP^{2}\setminus [1,0,0]$ with the normal bundle of $D_{1}$ and $\mathbb CP^{2}\setminus [0,1,0]$ with the normal bundle of $D_{2}$, these choices suffice to define $(\mathbb CP^{2},D_{1}\cup D_{2})$ as a manifold with a normal-crossing divisor in the sense of Definition 1.3. With this  choice, Ionel's notion of  rescaling by $\lambda$ towards $D_{1}$ corresponds to sending $[z_{1},z_{2},z_{3}]$ to $[\lambda^{-1} z_{1},z_{2},z_{3}]$, and rescaling by $\lambda$ towards $D_{2}$ corresponds to sending  $[z_{1},z_{2},z_{3}]$ to $[ z_{1},\lambda^{-1} z_{2},z_{3}]$.

\

Our first example has a rather complicated limit, and illustrates some  clarifications  necessary for the  understanding of \cite{IonelGW}.

Consider the the following sequence of maps  $f_{n}:\mathbb CP^{1}\longrightarrow \mathbb CP^{2}$
\[f_{n}([w_{1},w_{2}])=[n^{- 4}w_{1},n^{-3}(w_{1}+w_{2}),w_{2}]\]
Ionel's picture\footnote{In this paper, we draw these pieces as if they had a symplectic form compatible with the $(\mathbb C^{*})^{2}$--action, so the image of each piece is a moment map. The image of holomorphic curves under such a moment map is an `amoeba'  looking roughly like the red curve drawn.
} for the limit of this sequence is as follows:

\includegraphics{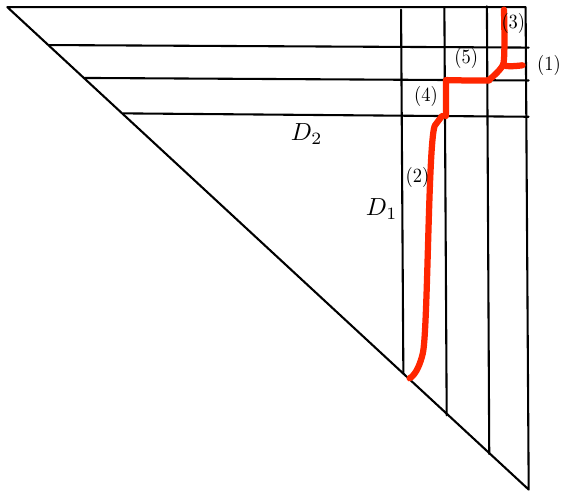}

There are $5$ interesting curve pieces in this picture: 
\begin{enumerate}
\item The curve in level $(3,2)$ is seen by rescaling the target by $n^{-4}$ towards $D_{1}$ and $n^{-3}$ towards $D_{2}$, 
\item The curve in level $(1,0)$ is seen by rescaling the target by $n^{-1}$ towards $D_{1}$ and rescaling the domain by $w_{1}\mapsto n^{3}w_{1}$.
\item The curve in level $(3,3)$ is a `trivial cylinder'  seen by rescaling the domain towards $[w_{1},w_{2}]=[-1,1]$ by a factor of $n^{-1}$, and rescaling the target by $n^{-3}$ towards both $D_{1}$ and $D_{2}$. 
\item The curve  between level $(1,1)$ and $(2,1)$ is a `trivial cylinder' seen by rescaling the domain by $w_{1}\mapsto n^{2} w_{1}$ and rescaling the target by $n^{-2}$ towards $D_{1}$ and  $n^{-1}$ towards $D_{2}$; to see this as a curve in the boundary of level $(1,1)$ we rescale instead by $n^{-1}$ towards $D_{1}$, and to see this as a curve in the boundary of level $(2,1)$, we rescale instead by $n^{-3}$ towards $D_{1}$. 
\item The curve between level $(2,1)$ and level $(2,2)$ is another `trivial cylinder' seen by rescaling the domain by $w_{1}\mapsto n w_{1}$, and rescaling the target by $(n^{-3},n^{-2})$ towards $(D_{1},D_{2})$; again, to see this as a curve in the boundary of level $(2,1)$ we rescale the target by $(n^{-3},n^{-1})$, and to see it as a curve in the boundary of level $(2,2)$, we rescale the target by $(n^{-3},n^{-3})$.
 \end{enumerate}
  To achieve the required rescaling of the target using Ionel's setup, we must rescale at each $n$ first by $n^{-1}$, then $n^{-2}$, then $n^{-1}$ so that what is seen in the limit on level $(a,b)$ corresponds to  rescaling  the target by $(n^{-\phi(a)},n^{-\phi(b)})$ where $\phi(0)=0,\ \phi(1)=1,\ \phi(2)=3, \ \phi(3)=4$.

Ionel's  Definition 7.3, of a limiting curve being relatively stable, must be interpreted (rather non-intuitively) using Definition 4.16  from \cite{IonelGW}.\footnote{See also the second paragraph of Remark 7.4. Be warned that the first paragraph of Remark 7.4 is  misleading. The correct condition for relative stability in the case of rescaling independently in different directions is the following: if one has to rescale $k$ times in the $i$th direction, then for each positive integer $0<m\leq k$, there is  a nontrivial curve in at least one multilevel that has $m$ in the $i$th spot. } In particular, the above limiting curve is relatively stable because there are nontrivial curves in levels $(1,0)$ and $(3,2)$, so, for the correct interpretation of Definition 7.3,  there are nontrivial curves in each positive level $1$, $2$, and $3$. Note that if we disliked curves being `between levels', we could avoid such in-between curves by rescaling $4$ times using a factor of $n^{-1}$ each time. The resulting curve would not be relatively stable. Furthermore, there exist simple examples  requiring infinitely many rescalings  to avoid all in-between cylinders.

   To satisfy Ionel's refined matching condition, (6.29) of  \cite{IonelGW},  between the curves in level $(1,0)$ and $(3,2)$, it is important that (at least in the limit),  the same rescaling factor is used the first and the third time: for example rescaling by $n^{-1}$ then $n^{-2}$, then $3n^{-1}$ would violate Ionel's refined matching condition (6.29), however rescaling by $2n^{-1}$ then $3n^{-2}$, then $2n^{-1}$ would be fine.  This essential point  is easily missed when intuitively considering the geometric pictures drawn in \cite{IonelGW}.

 Ionel's other refined matching condition --- the existence of a strictly negative simultaneous solution to the system of equations (6.25) --- is easily misinterpreted when there are in-between cylinders contained in the singular divisors between levels. We must  only apply (6.25) to positive-depth nodes or chains of trivial holomorphic cylinders, `$y$', connecting components with well-defined contact information $a_{i}$ in the $i$th direction\footnote{To interpret the statement of (6.25) correctly, the reader must be aware that for a node or chain of trivial holomorphic cylinders, `$y$', Ionel's  definition of $I^{\pm}(y)$  --- contained in Remark 5.12 --- includes only the directions in which this geometric contact information is well defined. A separate issue is that (6.25), as literally stated, also does not apply when $s_{i}=0$,  but rearranging to avoid dividing by $s_{i}$ gives a statement that obviously does apply. This trivial observation is needed for the correctly-interpreted set of equations from (6.25) to be equivalent to the equations for existence of certain tropical curves.} --- for example, in direction $1$ there is no equation using the node between components (2) and (4), but there is an equation using the nodes between components (2) and (5).\footnote{The astute reader looking at the proof of (6.25) may detect that, when contact information in the $i$th direction is not well defined, (6.25) is replaced by a pair of strict inequalities corresponding to at least one curve being in-between levels in the $i$th direction. There is no need to keep track of these inequalities because they follow from the equations for other nodes or chains of cylinders for which contact information is well defined.} In fact, a choice of strictly-negative simultanous solution to the equations (6.25) for all nodes or chains of trivial holomorphic cylinders `$y$' and all directions $i\in I^{\pm}(y)$ is equivalent to a choice of a certain tropical curve.  In our current example, following the notation of equation (6.25),  the system of equations is as follows: 

Assign variables
\begin{enumerate}
\item $\alpha(1)$ to the node between curve  (2) and (4),
\item $\alpha(2)$ to the node between curve (4) and (5),
\item $\alpha(3)$ to the node between curve (5) and (1),
\item $\alpha(4)$ to the node between curve (1) and (3)
\item $\alpha_{1}$ $\alpha_{2}$ and $\alpha_{3}$  to levels $1$, $2$ and $3$ respectively
\end{enumerate}

The equations from the direction $i=1$ are:
\[\alpha(1)+\alpha(2)=\alpha_{2}\]
\[\alpha(1)+\alpha(2)+\alpha(3)=\alpha_{2}+\alpha_{3}\]
\[\alpha(3)=\alpha_{3}\]
The equations from the direction $i=2$ are:
\[\alpha(1)=\alpha_{1}\]
\[\alpha(1)+\alpha(2)+\alpha(3)=\alpha_{1}+\alpha_{2}\]
\[\alpha(2)+\alpha(3)=\alpha_{2}\]
\[\alpha(4)=\alpha_{3}\]
There is an obvious two-dimensional family of strictly negative solutions: $\alpha(1)=\alpha(3)=\alpha(4)=\alpha_{1}=\alpha_{3}=a$,  $\alpha(2)=b$, and $\alpha_{2}=a+b$. The above are also  the equations for the existence of the tropical curve pictured below. Abstractly, this tropical curve is the dual graph of the nodal curve pictured above, however
\begin{itemize}\item it is given a complete metric so that the internal edges have lengths $-\alpha(i)$,
\item it is mapped to $[0,\infty)^{2}$ so that a vertex corresponding to a curve on level $(m,n)$ has first coordinate $-\sum_{i\leq m}\alpha_{i}$ and second  coordinate $-\sum_{i\leq n}\alpha_{i}$
\item and the derivative on each edge is an integral vector corresponding to Ionel's contact information $(\epsilon_{1}s_{1},\epsilon_{2}s_{2})$.\footnote{I am uncertain of this point, as I found it difficult to understand Ionel's exact definition of contact information $s_{i}$ with sign $\epsilon_{i}$, and this interpretation may be at odds with our above interpretation of what contact information is defined for the purposes of (6.25).  What is certain is that this integral vector is determined from Ionel's contact information by some such simple formula.}
\end{itemize}

\includegraphics{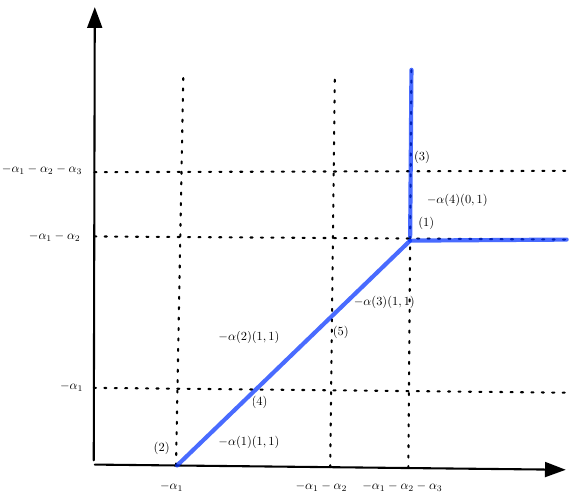}

This tropical curve is a Gromov--Hausdorf limit of the image of our sequence of curves $f_{n}$. Choose a metric $g$ on $\mathbb CP^{2}\setminus(D_{1}\cup D_{2})$ so that the real and imaginary parts of $z_{i}\frac \partial{\partial z_{i}} $ are orthonormal near $D_{i}$, and so that $g(v_{1},v_{2})$ is a smooth function for any pair of smooth vectorfields $v_{i}$ on $\mathbb CP^{2}$  tangent to both $D_{1}$ and $D_{2}$. Consider the image of $f_{n}$ in $\mathbb CP^{2}\setminus (D_{1}\cup D_{2})$ with the rescaled metric $(\log n)^{-1}g$. The Gromov--Hausdorf limit of this sequence of metric spaces is the quadrant pictured above: everything in a compact subset of $\mathbb CP^{2}\setminus(D_{1}\cup D_{2})$ is sent to the origin, traveling right corresponds to traveling towards $D_{1}$, and travelling up corresponds to travelling towards $D_{2}$. More specifically, given any point $[p_{1},p_{2},1]$ with $p_{i}\neq 0$, the sequence of points $[n^{-x}p_{1},n^{-y}p_{2},1]$ has limit $(x,y)$ above.  The blue tropical curve is the limit of the image of $f_{n}$.

There is an obvious $2$-dimensional family of tropical curves with the data required by our equations. We can parameterize the moduli space of these tropical curves by the location of the trivalent vertex, which must be located in the strictly positive span of $(1,1)$ and $(3,2)$ for the dotted lines above to intersect it as shown. The two dimensions of freedom of these tropical curves corresponds to the fact that we must quotient Ionel's moduli space of curves with this data by $(\mathbb C^{*})^{2}$, because all the different curves we identify in this way may be obtained by tweaking the rescaling parameters used while preserving Ionel's refined matching conditions. How each of our labeled vertices move under deforming the tropical curve corresponds to the effect of the $\mathbb C^{*}$--action on the piece of curve with the corresponding label. In particular, the curve (1) is rescaled by multiplying by $(c_{1},c_{2})$, the curve (2) is rescaled by multiplying by $(c_{1}/c_{2},1)$, and the curve (3) is rescaled by multiplying by $(c_{1},c_{1})$. After quotienting by this $(\mathbb C^{*})^{2}$--action, there is just one curve in Ionel's moduli space corresponding to the above pictures. 

The fact that that the vertex labeled by (1) may be located anywhere within the positive span of $(1,1)$ and $(3,2)$  corresponds to the  fact that the sequence of curves $f_{n}([w_{1},w_{2}]):=[x_{n}w_{1},y_{n}(w_{1}+w_{2}),w_{2}]$ has this limit if and only if $x_{n}\to 0$, $x_{n}/y_{n}\to 0$, and $y_{n}^{3}/x_{n}^{2}\to 0$ as $n\to 0$.  To see this, rescale by $x_{n}/y_{n}$, then $y^{2}_{n}/x_{n}$, then $x_{n}/y_{n}$. The position of the in-between trivial holomorphic cylinders is ensured by the condition that $y_{n}^{3}/x_{n}^{2}\to 0$.

It is easy to miss the importance of satisfying equation (6.25) or having the corresponding tropical curve. Curves violating this condition do not glue together, and are never the limits of holomorphic curves. Accordingly, tropical curves play a crucial   role in the gluing formulas for Gromov--Witten invariants that correspond to these relative invariants.

\

What is the advantage of thinking about the blue tropical curve instead of Ionel's equivalent equations (6.25)?  Given any sequence of lines that fall into $D_{1}\cap D_{2}$,  there is some sequence of rescalings resulting in a line in $(\mathbb C^{*})^{2}$ equal to  the curve (1) above. In the tropical picture, this curve will correspond to a trivalent vertex with edges leaving in directions $(1,0)$, $(0,1)$, and $(-1,-1)$.   A second (possibly trival) sequence of rescalings  will result in a monomial map of a cylinder into the relevant model for level $(a,0)$, $(0,b)$, or $(0,0)$ ---  corresponding tropically to a monovalent vertex with an edge leaving in direction $(1,1)$. For all other (successful)  rescalings,  the resulting curve will  be  a trivial holomorphic cylinder, corresponding tropically to a bivalent vertex with two  edges leaving  in opposite directions. We can ignore these for now. Solving for such tropical curves  is, in this case, a trivial task. Once we have a given tropical curve, we can fill in dotted lines as above,  then determine Ionel's pictures.

Solving for the tropical curve first, we  now draw representative pictures of all curves in Ionel's compactified moduli space. 

\

The following tropical picture

\includegraphics{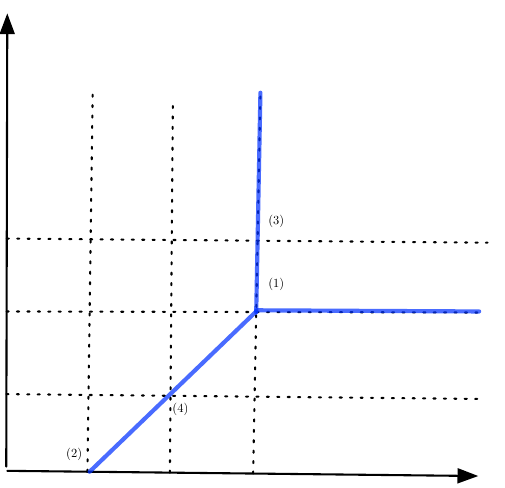}\\
corresponds to sequences
\[f_{n}([w_{1},w_{2}]):=[x_{n}w_{1},y_{n}(w_{1}+w_{2}),w_{2}]\]
with $x_{n}\to 0$ and $x_{n}^{2}/y_{n}^{3}\rightarrow a\neq 0$.

\

There is a one-parameter family of such tropical curves, so we quotient by $\mathbb C^{*}$ to obtain Ionel's moduli space --- the equations (6.25) for the tropical curves are obtained from the  equations determining the refined matching conditions;  each real freedom we have in solving (6.25) corresponds to a $\mathbb C^{*}$--worth of freedom in  matching-condition-compatible rescaling parameters. In practice, what quotient is needed is very easy to work out from the the tropical picture.  In this case, the action on curve pieces is multiplication of $(1)$ by  
 $(c^{3},c^{2})$, $(2)$ by $(c,1)$, 
$(3)$ by $(c^{3},c^{3})$,
  and $(4)$ by $(c^{2},c)$.

In this situation, there is a two-complex-parameter family of Ionel's curves, parametrized by the possible lines in $(\mathbb C^{*})^{2}$ corresponding to (1). After taking the quotient, Ionel's compactified modulis space has a one-complex-parameter family of such limit curves parameterized by the limit of $x_{n}^{2}/y_{n}^{3}$. Ionel's pictures of these limiting curves look as follows: 

\

\includegraphics{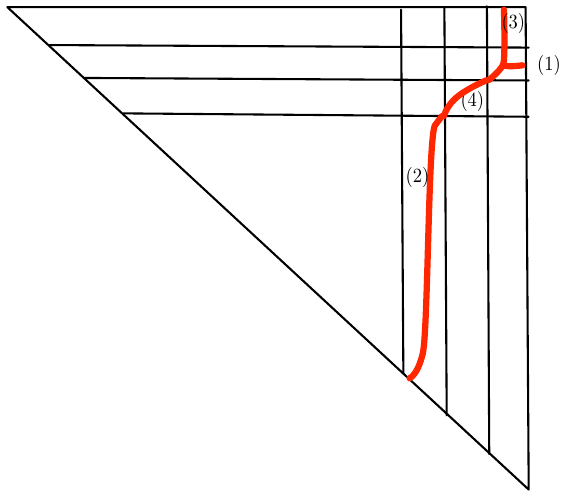}

To verify the above results using Ionel's method, rescale by $x_{n}/y_{n}$ three times.

\

The following tropical picture

\includegraphics{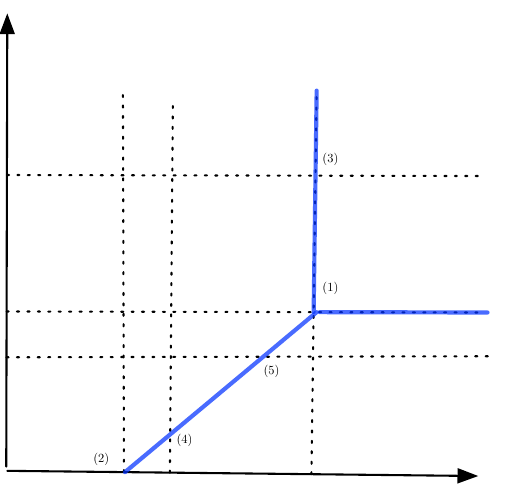}

 \noindent corresponds to  sequences 
$f_{n}([w_{1},w_{2}]):=[x_{n}w_{1},y_{n}(w_{1}+w_{2}),w_{2}]$
with $x_{n}\to 0$, $x_{n}^{2}/y_{n}^{3}\rightarrow 0$, and $y_{n}^{2}/x_{n}\rightarrow 0$.

There is a $2$--dimensional family of such tropical curves, and we must quotient by a $(\mathbb C^{*})^{2}$--action, multiplying $(1)$ by $(c_{1},c_{2})$, 
$(2)$ by $(c_{1}/c_{2},1)$,
 and $(3)$ by $(c_{1},c_{1})$.

Ionel's picture of such a limiting curve is below. After taking the quotient, there is only one such curve. 

\includegraphics{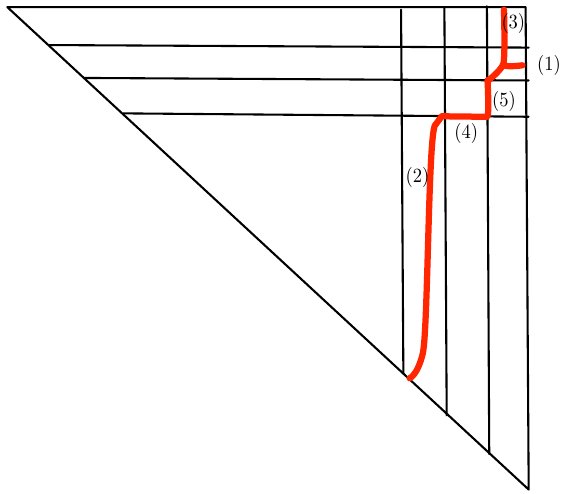}

To see this, rescale by $x_{n}/y_{n}$, $y_{n}^{2}/x_{n}$, then $x_{n}/y_{n}$.

\

For sequences 
$f_{n}([w_{1},w_{2}]):=[x_{n}w_{1},y_{n}(w_{1}+w_{2}),w_{2}]$
with $x_{n}\to 0$ and $x_{n}/y_{n}^{2}\rightarrow a\neq 0$, we get the following tropical picture.

\includegraphics{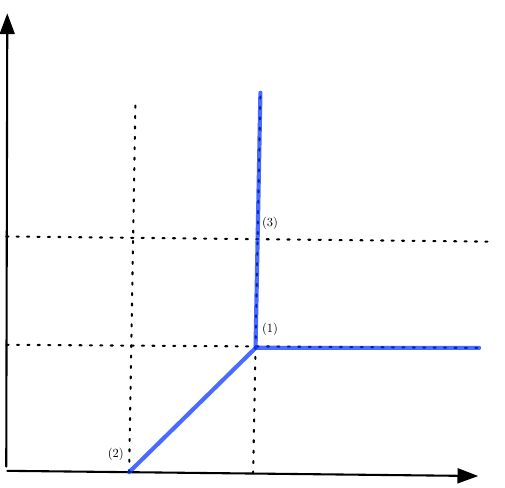}

There is a $1$--dimensional family of such tropical curves, and we quotient by a $\mathbb C^{*}$--action multiplying
$(1)$ by $(c^{2},c)$, 
$(2)$ by $(c,1)$,
and $(3)$ by $(c^{2},c^{2})$.

After taking this quotient, there is a $1$-complex dimensional family of such limiting curves parameterized by the limit of $x_{n}/y_{n}^{2}$. Ionel's pictures of such curves are as follows:

\includegraphics{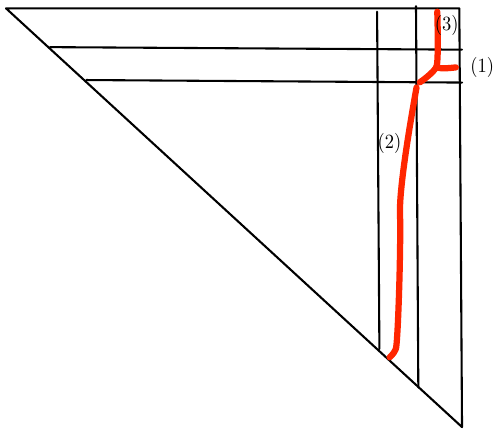}

To see this, rescale twice by $y_{n}$.

\

For sequences 
$f_{n}([w_{1},w_{2}]):=[x_{n}w_{1},y_{n}(w_{1}+w_{2}),w_{2}]$
with $y_{n}\to 0$ and $x_{n}/y_{n}^{2}\to 0$, we get the following tropical picture.

\includegraphics{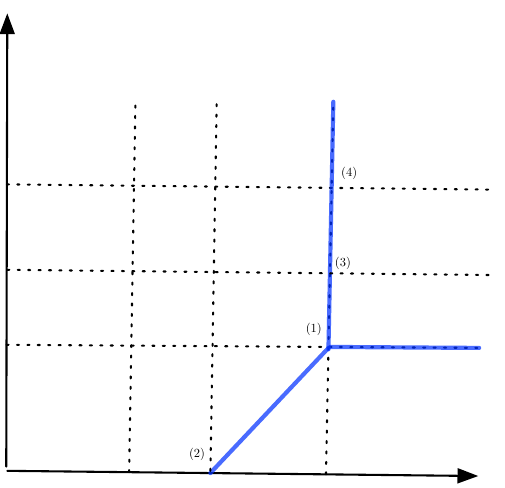}

There is a $2$--dimensional family of such tropical curves, and we quotient by $(\mathbb C^{*})^{2}$, to obtain a unique curve in Ionel's moduli space with this data. Ionel's picture of such a curve is as follows:

\includegraphics{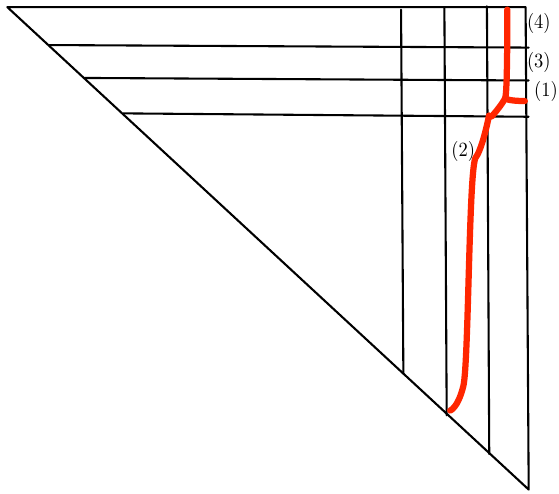}

To see this, rescale by $y_{n}$, then $x_{n}/y_{n}^{2}$, then $y_{n}$.

\

For sequences 
$f_{n}([w_{1},w_{2}]):=[x_{n}w_{1},y_{n}(w_{1}+w_{2}),w_{2}]$
with $x_{n}\to 0$ and $y_{n}\to a\neq 0$, we get the following tropical picture.

\includegraphics{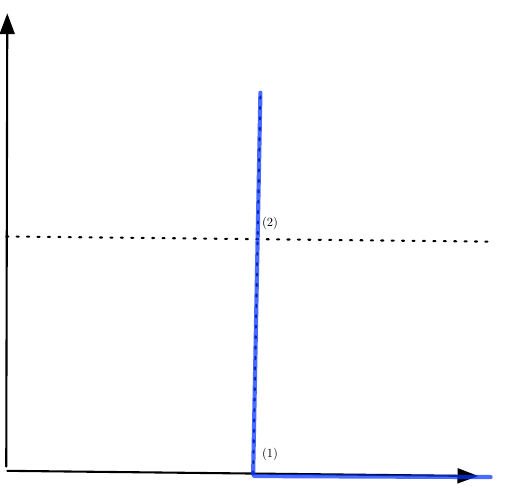}

There is a $1$--dimensional family of such tropical curves, so we quotient by $\mathbb C^{*}$ to obtain a 
 $1$--complex-dimensional family of limiting curves in Ionel's moduli space, parameterised by the limit of $y_{n}$. Ionel's picture for such curves is as follows: 

\includegraphics{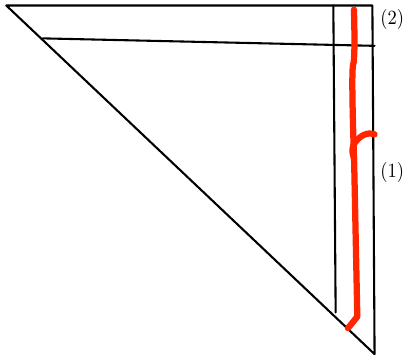}

All the above limiting curves have mirror-image versions  swapping the role of $x_{n}$ and $y_{n}$. There are also the following curves where $x_{n}$ and $y_{n}$ play a roughly symmetric role:

For sequences  $f_{n}([w_{1},w_{2}]):=[x_{n}w_{1},y_{n}(w_{1}+w_{2}),w_{2}]$
with $x_{n}\to 0$ and $x_{n}/y_{n}\to a\neq 0$, we get the following tropical picture.

 \includegraphics{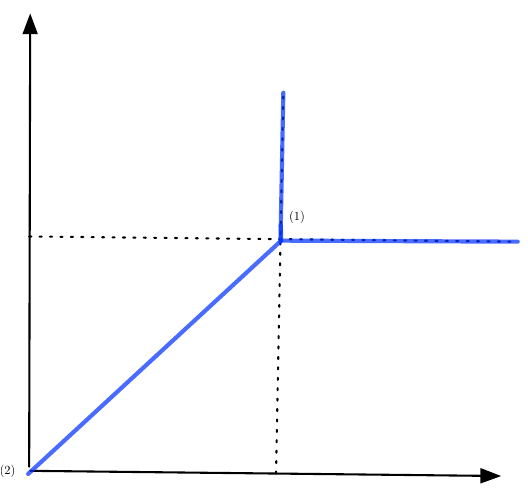}
 
 There is a $1$--dimensional family of such tropical curves, so we  quotient by $\mathbb C^{*}$. After this quotient, Ionel's moduli space contains a $1$-complex dimensional family of such curves, parameterised by the limit of $x_{n}/y_{n}$.
Ionel's picture of such curves is as follows:

\includegraphics{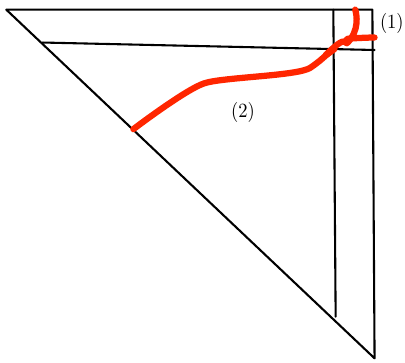}

\

Before continuing with Ionel's moduli space, we quickly describe what moduli space is obtained using exploded manifolds.

$\mathbb CP^{2}$ with the normal-crossing divisor $D_{1}\cup D_{2}$ is a complex manifold with a normal-crossing divisor, so the explosion functor from \cite{iec} applies, giving an exploded manifold $\ex B$. Like all exploded manifolds, $\ex B$ has a tropical part, $\totb{\ex B}=[0,\infty)^{2}$, spookily similar to the quadrant  in the above  tropical pictures.  Each holomorphic curve in $\ex B$ has a tropical part --- a tropical curve in $\totb{\ex B}$. The tropical parts of the relevant curves in $\ex B$  are the blue tropical curves drawn above, however now we do not care about the dotted lines corresponding to Ionel's levels. The moduli space of such tropical curves is parameterised by the position of the central vertex. Because the tropical curves of the last type drawn above are special, we divide this moduli space along the corresponding ray $(1,1)$ so that the moduli space of these tropical curves is the following fan:

\includegraphics{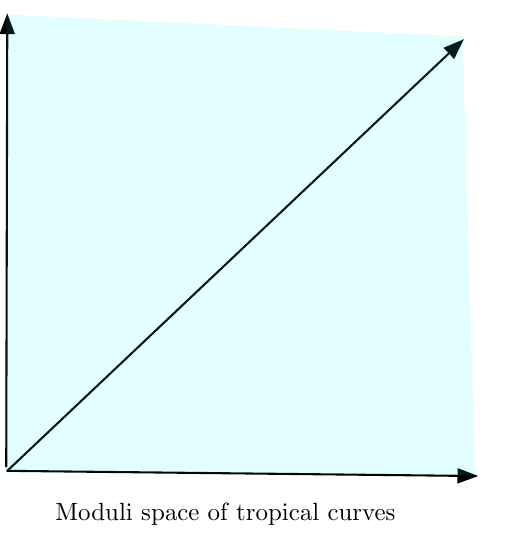}

The moduli space of relevant curves in $\ex B$ is again an exploded manifold $\ex M$. The tropical part $\totb{\ex M}$ of this moduli space is the above moduli space of tropical curves.\footnote{In more general situations, the tropical part of the moduli space of curves is some union of moduli spaces of tropical curves, so things are a little more complicated.} In this case, the moduli space $\ex M$ itself is easily identified --- it is the explosion of the toric manifold with the following moment map relative to the normal-crossing divisor made up of the three upper-right boundary divisors:

\

\includegraphics{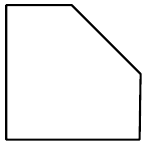}

\

The above toric manifold may either be viewed as the moduli space of lines in $\mathbb CP^{2}$ blown up at the points corresponding $D_{1}$ and $D_{2}$, or as the blow up of $D_{1}\times D_{2}$ at the point corresponding to their intersection. From the first perspective, the intersection point of $D_{1}$ and $D_{2}$ gives a dual line in the moduli space of lines, and the normal-crossing divisor is the inverse image of this dual line under the blowup. From the second perspective, the normal-crossing divisor is the inverse image of $(\{0\}\times D_{2})\cup (D_{1}\times\{0\})$ under the blowup map.

\

What about Ionel's compactified moduli space? If we pay attention to the dotted lines in the tropical pictures above,  we must subdivide our moduli space of tropical curves as follows:

\includegraphics{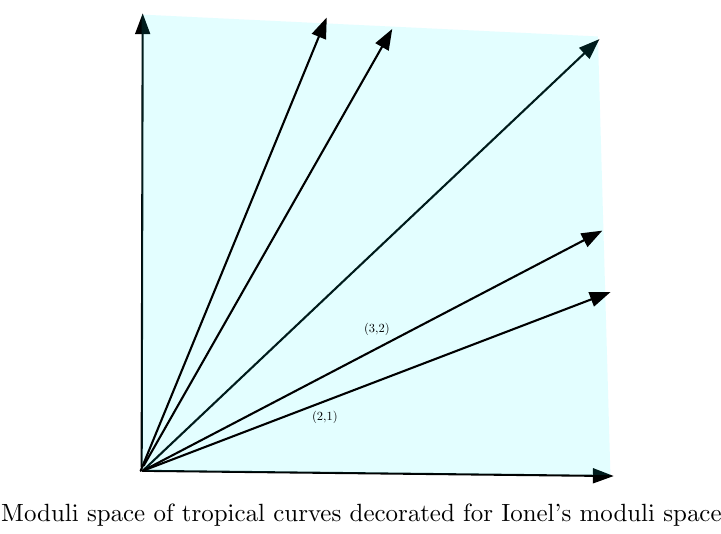}

 Ionel's compactified moduli space is the toric manifold with toric fan the union of the above fan with the rays going down and left; the extra strata added in during the compactification are the boundary strata corresponding to the toric fan above. To get Ionel's moduli space from the toric manifold with moment map pictured above,  blow up the two upper right corners, then again blow up the two upper right corners of the resulting toric manifold. The last two blowups could be avoided by  using a version of Ionel's moduli space  ignoring the (largely irrelevant) information of where trivial holomorphic cylinders go. All four blowups could be avoided by using a version of Ionel's moduli space ignoring Ionel's `level' labelling.

Ionel's enhanced evaluation map in this case is the blowdown map to $D_{1}\times D_{2}$, so Ionel's $4$ extra blowups  are irrelevant to the enhanced evaluation map. To perturb the moduli space of curves as allowed by Ionel,  we must modify $J$ compatibly with the normal-crossing divisors, or choose a smooth complex Donaldson--divisor $N$  intersecting $D_{1}\cup D_{2}$ transversely, then take note of where a curve intersects $N$, and perturb the right hand side of the $\dbar$ equation using this information.  As such perturbations take no notice of where trivial holomorphic cylinders are, or the extra `level' information in Ionel's moduli space, the enhanced evaluation map from the perturbed moduli space will still factor through the un-blown-up  moduli space.  It is not difficult to prove that the same is true in general.

\

We  now sketch the general  relationship between Ionel's compactified moduli space and the moduli spaces constructed using exploded manifolds.

 Each exploded manifold $\ex B$ has a smooth part, $\totl{\ex B}$, a toplogical space with smooth manifold strata (and a bit of extra `smooth' structure near strata.) For example, if $\ex B$ is the explosion of a complex manifold with normal-crossing divisors, $\totl{\ex B}$ is that complex manifold with normal-crossing divisors.  When $\totl{\ex B}$ can be given Ionel's structure of a symplectic manifold with normal-crossing divisors, we may compare moduli spaces. (This is the case if and only if $\ex B$ admits a symplectic structure and the tropical part of $\ex B$ is  a union of quadrants.) Starting from one of Ionel's symplectic manifolds with normal-crossing divisors, we must make sure that the almost complex structure $J$ used satisfies a mild assumption called $\dbar \log$--compatibility  (described in section 14 of \cite{elc}). Then, there exists an exploded manifold $\ex B$ whose  smooth part $\totl{\ex B}$ is this symplectic manifold with normal-crossing divisors.

Ionel's Definition 7.2 of a map into a level $m$ building is a naive map (Definition 7.1) which also satisfies the refined matching condition (6.29), and which has a balanced\footnote{Tropical geometers be warned: this is Ionel's definition of balanced, it is unrelated to the balancing condition that tropical curves satisfy at vertices. Also be warned that the tropical curves coming from exploded manifolds or log geometry are not necessarily balanced in the tropical geometer's sense.} refined dual graph --- the existence of such a dual graph amounts to the existence of strictly negative solutions to the set of equations (6.25), correctly interpreted as above. A low tech way of specifying any curve in the corresponding exploded manifold to specify such a map, choose a strictly negative solution to (6.25), choose a negative number\footnote{These extra negative numbers may be incorporated into the system of equations (6.25) by extending its applicability to all nodes, and not just positive depth nodes as in \cite{IonelGW}.} and number in $\mathbb C^{*}$  for each node with contact information $\{s_{i}\}$  the $0$--vector,  then choose  a little extra gluing information ---   a $k$-fold choice of gluing information for each node (or chain of trivial cylinders) with contact information $(s_{1},\dotsc,s_{n})$   $k$ times a primitive integer vector.  (The proof of this observation is trivial after the reader has understood the definition of exploded manifolds and the explosion functor, found in sections 3 and 5 of \cite{iec} --- the tropical part of the exploded curve is specified by the solution to (6.25) and the choice of negative number for each zero depth node, then  the remaining information specifies an exploded curve. The one complication is that  exploded curves have an extra regularity requirement at nodes, but this regularity is automatically satisfied by holomorphic curves.)

Conversely, we can construct any of Ionel's maps into a level $m$ building starting with exploded curves. The most effective method uses the notion of a refinement of an exploded manifold, described in section 10 of \cite{iec}. A refinement of an exploded manifold $\ex B$ may be specified by a subdivision of its tropical part, but the effect on the smooth part of $\ex B$ is a kind of blowup. For example, the moduli space of exploded curves $\ex M$ described above has  tropical part pictured above as the moduli space of tropical curves. $\ex M$  has a refinement corresponding to the subdivision pictured as the moduli space of tropical curves decorated for Ionel's moduli space. The smooth part of this refinement of $\ex M$ is Ionel's moduli space, which, as we have described, is a blowup of the smooth part of $\ex M$. Refinement is almost an isomorphism in the category of exploded manifolds. For example, it is proved in \cite{egw} that refining an exploded manifold doesn't really change Gromov--Witten invariants, because it just results in a refinement of the moduli space of curves. (Analogous results are proved in the log setting in \cite{ilgw}.)

 With refinements understood, the following procedure  constructs one of Ionel's maps from an exploded curve $f$ in $\ex B$. Every vertex of the tropical part $\totb f$ of $f$ has coordinates $(x_{1},\dotsc,x_{n})$ within some strata of $\totb{\ex B}$ isomorphic to a quadrant $[0,\infty)^{n}$. Let $0<l_{1}<\dotsb<l_{m}$ be the list  of nonzero numbers that occur in such coordinates. There is a refinement $\ex B'$ of $\ex B$ specified by subdividing the tropical part of $\ex B$ by all planes that have an $l_{i}$ as a coordinate, and  there is a corresponding refinement $f'$ of $f$, a refinement that is an exploded curve in $\ex B'$. Then,  the smooth part of $\ex B'$ is one of Ionel's level $m$ buildings and the smooth part of $f'$ is one of Ionel's maps into a level $m$ building. If $f$ is stable,  resulting map is relatively stable in the sense of Ionel's Definition 7.3 (and would not be relatively stable if we refined using any extra numbers in addition to  \{$l_{i}$\}).

 For example, the dotted lines in the tropical pictures above determine a refinement with smooth part Ionel's corresponding building. For each of these pictures, there is an exploded curve $f$ whose tropical part is the blue tropical curve depicted. The image of the smooth part of $f'$ is the corresponding red curve depicted in Ionel's pictures.

\

The relationship between the moduli space of exploded curves and Ionel's moduli space is also most effectively described using refinements. The moduli space of tropical curves in $\totb{\ex B}$ can be subdivided into strata  depending on how many of the coordinates of their vertices coincide. An example is the subdivision above of our moduli space of tropical curves that give different pictures for Ionel's curves. This subdivision of the moduli space of tropical curves defines a refinement $\ex M'$ of the moduli space of curves, $\ex M$. Ionel's moduli space is the smooth part  $\totl{\ex M'}$ of $\ex M'$, and is therefore  a blow up of the smooth part $\totl{\ex M}$ of the moduli space of exploded curves.\footnote{There is one technical issue with this assertion that I am not totally clear on: When extra gluing information is required to describe the exploded moduli space because of trivial cylinders that are multiple covers, there may be some kind of orbifold information that is not translated correctly when passing from the exploded moduli space to Ionel's moduli space. I think that this problem is taken care of by Ionel taking a quotient using the torus $T_{\sigma}$ instead of the torus $\tilde T_{\sigma}$ described above Definition 7.2.}

\bibliographystyle{plain}
\bibliography{ref.bib}
\end{document}